\documentclass[11pt, a4paper]{amsart}
\usepackage{mathbbol}
\usepackage{amsmath, amsthm, amssymb}
\usepackage{txfonts, mathrsfs}
\usepackage[usenames]{color}
\usepackage{psfrag}
\usepackage{graphicx}

\newtheorem*{thm*}{Theorem}

\newtheorem{thm}{Theorem}[section]
\newcommand{\bt}{\begin{thm}}
\newcommand{\et}{\end{thm}}

\newtheorem{cor}[thm]{Corollary}
\newcommand{\bc}{\begin{cor}}
\newcommand{\ec}{\end{cor}}

\newtheorem{lem}[thm]{Lemma}
\newcommand{\bl}{\begin{lem}}
\newcommand{\el}{\end{lem}}

\newtheorem{prop}[thm]{Proposition}
\newcommand{\bp}{\begin{prop}}
\newcommand{\ep}{\end{prop}}

\newtheorem{defn}[thm]{Definition}
\newcommand{\bd}{\begin{defn}}      
\newcommand{\ed}{\end{defn}}

\newtheorem{rmrk}[thm]{Remark}
\newcommand{\br}{\begin{rmrk}}
\newcommand{\er}{\end{rmrk}}

\newtheorem{quest}[thm]{Question}
\newcommand{\bq}{\begin{quest}}
\newcommand{\eq}{\end{quest}}

\newcommand{\N}{\mathbb{N}}

\newcommand{\R}{\mathbb{R}}

\newcommand{\Hyp}{\mathbb{H}}

\newdimen\vintkern\vintkern12pt
\def\vint{-\kern-\vintkern\int}

\newcommand{\trace}{\operatorname{tr}}

\newcommand{\hm}{{\mathcal H}}

\begin{document}
\bibliographystyle{plain}

\title[Harmonic quasi-isometric maps]{Harmonic quasi-isometric maps into Gromov hyperbolic ${\rm CAT}(0)$-spaces}

\author{Hubert Sidler}

\address
  {Department of Mathematics\\ University of Fribourg\\ Chemin du Mus\'ee 23\\ 1700 Fribourg, Switzerland}
\email{hubert.sidler@unifr.ch}

\author{Stefan Wenger}

\address
  {Department of Mathematics\\ University of Fribourg\\ Chemin du Mus\'ee 23\\ 1700 Fribourg, Switzerland}
\email{stefan.wenger@unifr.ch}

\thanks{Research partially funded by Swiss National Science Foundation Grant 165848.}

\date{\today}

\begin{abstract}
 We show that for every quasi-isometric map from a Hadamard manifold of pinched negative curvature to a locally compact, Gromov hyperbolic, ${\rm CAT}(0)$-space there exists an energy minimizing harmonic map at finite distance. This harmonic map is moreover Lipschitz. This generalizes a recent result of Benoist-Hulin.
\end{abstract}

\maketitle

\section{Introduction}

The well-known Schoen-Li-Wang conjecture asserts that every quasiconformal self-homeomorphism of the boundary at infinity of a rank one symmetric space $M$ extends to a unique harmonic map from $M$ to itself. This conjecture has recently been settled in the affirmative in a series of break-through papers by Markovic \cite{Mar15}, \cite{Mar17}, Lemm-Markovic \cite{LM18}, and Benoist-Hulin \cite{BH-1-17}.  Earlier partial results were proved in \cite{Pan89}, \cite{TW95}, \cite{HW97}, \cite{Mar02}, \cite{BS10}, see also the references in \cite{BH-1-17}.
Benoist-Hulin's result \cite{BH-1-17}, which goes beyond the Schoen-Li-Wang conjecture, shows that every quasi-isometric map between rank one symmetric spaces $X$ and $Y$ is at finite distance of a unique harmonic map. 
Even more recently, Benoist-Hulin  \cite{BH-2-17} extended their result in \cite{BH-1-17} to the case when $X$ and $Y$ are Hadamard manifolds of pinched negative curvature, i.e.~simply connected Riemannian manifolds of sectional curvature bounded by $-b^2\leq K_X, K_Y \leq -a^2$ for some constants $a,b>0$. 

The aim of the present note is to further generalize the existence part of Benoist-Hulin's result  \cite{BH-2-17} by relaxing the curvature conditions on the target space $Y$. Our methods even work in the context of singular metric spaces $Y$. 
Recall that Korevaar-Schoen \cite{KS93} developed a theory of Sobolev and harmonic maps from a Riemannian domain into a complete metric space. We refer to \cite{KS93} and to Section~\ref{sec:Sobolev} of the present note for the definition. Our main theorem is:
\bt\label{thm:main-intro}
 Let $X$ be a Hadamard manifold of pinched negative curvature and let $Y$ be a locally compact, Gromov hyperbolic, $\mathrm{ CAT}(0)$-space. Then for every quasi-isometric map $f\colon X\to Y$ there exists an energy minimizing harmonic map $u\colon X\to Y$   which is globally Lipschitz continuous and at bounded distance from $f$. 
 \et

The precise meaning of being at bounded distance from $f$ is that $$\sup_{x\in X} d_Y(u(x), f(x))<\infty.$$ It follows in particular that $u$ is also quasi-isometric. Recall that a map $f\colon X\to Y$ between metric spaces $(X, d_X)$ and $(Y, d_Y)$ is called $(L,c)$-quasi-isometric if $$L^{-1} \cdot d_X(x,x') - c \leq d_Y(f(x), f(x')) \leq L\cdot d_X(x,x') + c$$ for all $x,x'\in X$. The map $f$ is called quasi-isometric if it is $(L,c)$-quasi-isometric for some $L\geq 1$ and $c\geq 0$. A quasi-isometric map is thus biLipschitz at large scales but no restriction is posed on small scales. In particular, $f$ need not be continuous. Notice moreover that the image $f(X)$ need not be quasi-dense in $Y$.

Recall that a geodesic metric space $Y$ is called ${\rm CAT}(0)$ if geodesic triangles in $Y$ are at least as thin as their Euclidean comparison triangles. Every Hadamard manifold is ${\rm CAT}(0)$. A geodesic metric space $Y$ is called Gromov hyperbolic if there exists $\delta\geq0$ such that each side of a geodesic triangle in $Y$ lies in the $\delta$-neighborhood of the other two sides. This is a large scale notion of negative curvature. It poses no restriction on small scales. We refer for example to \cite{BH99}, \cite{Gro87}, \cite{GdH90} for comprehensive accounts on ${\rm CAT}(0)$-spaces and Gromov hyperbolicity. Since every Hadamard manifold $Y$ of pinched negative curvature is locally compact, Gromov hyperbolic, and ${\rm CAT}(0)$ our Theorem~\ref{thm:main-intro} in particular recovers the existence part of Benoist-Hulin's result \cite[Theorem 1.1]{BH-2-17}. 

Unlike in the setting of \cite{BH-2-17}, energy minimizing harmonic maps at finite distance from a fixed quasi-isometric map need not be unique in our more general setting. Indeed, if $X=\Hyp^2$ is the hyperbolic plane and $Y:= \Hyp^2\times [0,1]$ then the maps $u_t(z):= (z,t)$ for $t\in [0,1]$ are isometric and energy minimizing harmonic and have finite distance from each other. In the context of singular target spaces uniqueness was shown in \cite{LW98} for harmonic maps at finite distance from a quasi-isometry between a cocompact Hadamard manifold and a cocompact ${\rm CAT}(\kappa)$-space with $\kappa<0$.

The main strategy of proof of our Theorem~\ref{thm:main-intro} is the same as that in \cite{BH-2-17}, and many of our arguments are in fact similar to those in \cite{BH-2-17}. On the one hand, existence and (local) Lipschitz regularity of energy minimizing harmonic maps is known in our more general context, see \cite{KS93}. On the other hand, the smooth structure of the target space $Y$ and the pinched negative curvature condition on $Y$ are crucially used at several places in \cite{BH-2-17}. This is for example essential when establishing bounds on the distance between a quasi-isometric map $f$ and a harmonic map. One of the principal new ingredients in our proof of similar bounds in our more general context is the use of the  Bonk-Schramm embedding theorem \cite{BS00}. This together with an argument about injective hulls, essentially due to \cite{Lan13}, allows us to rough-isometrically embed the (non-geodesic) image $f(X)$ into the hyperbolic $k$-space $\Hyp^k$ of constant curvature $-1$ for some $k\in\N$. The rough-isometric condition, which is much stronger than the quasi-isometric condition, then allows us to prove estimates on the distance between a quasi-isometric map $f$ and an energy minimizing harmonic map similarly to \cite{BH-2-17}. A further but more minor difference between our arguments and those in \cite{BH-2-17} is that we consistently work with the Gromov product in the target space $Y$ whereas the arguments in \cite{BH-2-17} rely on an interplay between estimates on the Gromov product and angle estimates. Such estimates on angles are not available in our setting since they require a strictly negative upper curvature bound.

\section{Preliminaries}\label{sec:preliminaries}

\subsection{Basic notation}

All metric spaces in our text will be complete. Let $(X, d)$ be a metric space. The open and closed balls in $X$ centered at $x\in X$ and of radius $r>0$ are denoted by $B(x,r):= \{x'\in X: d(x,x')<r\}$ and $\bar{B}(x, r):= \{x'\in X: d(x, x')\leq r\}$, respectively. The distance sphere is $S(x,r):= \{x'\in X: d(x,x')= r\}$. The Hausdorff $n$-measure on $X$ will be denoted by $\hm^n$. The normalization factor is chosen in such a way that $\hm^n$ equals the Lebesgue measure on Euclidean $\R^n$. In particular, if $X$ is a Riemannian manifold of dimension $n$ then $\hm^n$ equals the Riemannian volume. The averaged integral will be denoted by $$\vint_Af\,d\hm^n:= (\hm^n(A))^{-1}\cdot\int_Af\,d\hm^n.$$

\subsection{Some Riemannian preliminaries}\label{sec:Riem-prelim}

Let $M$ be a Riemannian manifold. The differential of a smooth function $f\colon M\to \R$ will be denoted by $Df$. The hessian $D^2f$ of $f$ is the $2$-tensor satisfying $$D^2f(X, X') = X(X'(f)) - (\nabla_X X')(f)$$ for all vector fields $X$, $X'$ on $M$. The trace of the hessian of $f$ is the Laplace of $f$ and denoted $\Delta f$. The function $f$ is called harmonic if $\Delta f\equiv 0$. If $M$ is an $n$-dimensional Hadamard manifold of sectional curvature $-b^2\leq K_M\leq -a^2$ for some $a,b>0$ then the hessian of the distance function $d_{x_0}$ to a point $x_0\in M$ satisfies 
\begin{equation}\label{eq:hessian-distance-hyp}
a \coth(a d_{x_0}) \cdot (g - D d_{x_0} \otimes D d_{x_0}) \leq D^2 d_{x_0}\leq b \coth(b d_{x_0}) \cdot (g - D d_{x_0} \otimes D d_{x_0})
\end{equation}
 on $M\setminus\{x_0\}$, where $g$ denotes the Riemannian metric on $M$. This follows from the hyperbolic law of cosines and comparison estimates, see e.g.~\cite{BH-1-17}. In particular, the laplacian of $d_{x_0}$ satisfies $\Delta d_{x_0}\geq a\cdot (n-1)$ on $M\setminus\{x_0\}$. 

Let $\varphi\colon M\to N$ be a smooth map into another Riemannian manifolds $N$. We denote by $D\varphi$ the differential of $\varphi$. The second covariant derivative of $\varphi$ is the vector-valued $2$-tensor which satisfies $$D^2\varphi (X, X') = \overline{\nabla}_X(D\varphi(X')) - D\varphi(\nabla_X X')$$ for all vector fields $X, X'$ on $M$, where $\overline{\nabla}$ denotes the pullback under $\varphi$ of the Riemannian connection on $N$. The trace of $D^2\varphi$ is called the tension field of $\varphi$ and denoted $\tau(\varphi)$.  If $\varphi\colon M\to N$ and $h\colon N\to \R$ are smooth then one calculates that 
\begin{equation}\label{eq:laplace-composition}
 \Delta(h\circ \varphi) = Dh(\tau(\varphi)) + \sum_{i=1}^n D^2h(D\varphi(e_i), D\varphi(e_i)),
\end{equation}
 where $\{e_1,\dots, e_n\}$ is an orthonormal basis in a tangent space of $M$.
The map $\varphi$ is called harmonic if $\tau(\varphi) \equiv0$.

\subsection{Gromov hyperbolicity}\label{sec:Gromov-hyp}

Let $(Y, d)$ be a metric space. Recall that the Gromov product of $x, y\in Y$ with respect to a basepoint $w\in Y$ is defined by 
\begin{equation}\label{eq:Gromov-product}
 (x\mid y)_w:= \frac{1}{2}\left[d(x,w) + d(y, w) - d(x,y) \right].
\end{equation}

\bd\label{def:Gromov-hyp}
A metric space $Y$ is called $\delta$-hyperbolic, $\delta\geq 0$, if $$(x\mid z)_w \geq \min\left\{(x\mid y)_w, (y\mid z)_w\right\} - \delta$$ for all $x,y,z,w\in Y$. The space is called Gromov hyperbolic if it is $\delta$-hyperbolic for some $\delta\geq 0$. 
\ed

A geodesic metric space $Y$ is Gromov hyperbolic in the sense of the definition above if and only if there exists $\bar{\delta}\geq 0$ such that every side of a geodesic triangle in $Y$ is contained in the $\bar{\delta}$-neighborhood of the other two sides, see \cite[Proposition III.H.1.22]{BH99}. 

For a proof of the following lemma see for example \cite[Theorem 3.21]{Vai05}.

\bl\label{lem:effecxt-qi-Gromov-product}
 Let $f\colon X\to Y$ be an $(L, c)$-quasi-isometric map between geodesic $\delta$-hyperbolic metric spaces $X$ and $Y$. Then there exists a constant $c'$ only depending on $L$, $c$, and $\delta$ such that for all $x,x', w\in X$ we have $$L^{-1}\cdot (x\mid x')_w - c'\leq (f(x)\mid f(x'))_{f(w)} \leq L\cdot (x\mid x')_w + c'.$$
\el

The next lemma is also known as exponential divergence of geodesics.

\bl\label{lem:div-geodesics}
 Let $(Y, d)$ be a geodesic $\delta$-hyperbolic metric space. Then there exists $\delta'>0$ depending only on $\delta$ with the following property. Let $r_1, r_2\geq 1$ and let $\gamma, \eta \colon [0, r_1+r_2]\to Y$ be two geodesics parametrized by arc-length with $\gamma(0) = \eta(0)$. If $d(\gamma(r_1), \eta(r_1))>3\delta'$ then any curve connecting $\gamma(r_1+r_2)$ and $\eta(r_1+r_2)$ outside the ball $B(\gamma(0), r_1+r_2)$ has length at least $2^{(r_2-1)\cdot \delta'^{-1}}$.
\el

This follows for example from the proof of \cite[Proposition III.H.1.25]{BH99}.

\subsection{Injective hulls of metric spaces}\label{sec:injective-hull}

We will need the following construction of an injective hull due to Isbell \cite{Isb64}. Given a metric space $(Z, d)$, denote by $E(Z)$ the space of all functions $f\colon Z\to\R$ satisfying 
\begin{equation}\label{eq:injective-hull}
 f(z) + f(z') \geq d(z,z')
\end{equation}
for all $z,z'\in Z$ and such that $f$ is extremal in the following sense. If $g\colon X\to\R$ is another function satisfying \eqref{eq:injective-hull} and $g\leq f$ then $g=f$. The space $E(Z)$, when equipped with the supremum norm, is called the injective hull of $Z$. It is an injective metric space in the sense that for every subset $A$ of a metric space $B$ and every $1$-Lipschitz map $\varphi\colon A\to E(Z)$ there exists a $1$-Lipschitz extension $\overline{\varphi}\colon B\to E(Z)$ of $\varphi$. In particular, it follows that $E(Z)$ is a geodesic metric space. The space $Z$ embeds isometrically into $E(Z)$ via the map $z\mapsto d(z,\cdot)$. Moreover, if $Z$ is a subset of another metric space $Z'$ then there exists an isometric embedding $h\colon E(Z)\to E(Z')$ such that $h(f)|_Z = f$ for all $f\in E(Z)$, see \cite[Proposition 3.5]{Lan13}. It was proved in \cite[Proposition 1.3]{Lan13} that if $Z$ is Gromov hyperbolic then so is $E(Z)$ and that if $Z$ is moreover geodesic then $E(Z)$ also lies in finite Hausdorff distance of $Z$.

\section{Sobolev maps into metric spaces}\label{sec:Sobolev}

There are several equivalent definitions of Sobolev maps from a Riemannian domain to a complete metric space, see for example \cite{HKST15} and the approaches described therein. We will use the definition given by Korevaar-Schoen in \cite{KS93}. As we will only deal with Sobolev maps of exponent $p=2$ and defined on open balls in a Hadamard manifold, we will restrict to this setting.

Let $X$ be a Hadamard manifold of dimension $n\geq 2$ and let $\Omega\subset X$ be an open, bounded ball. Let $(Y, d_Y)$ be a complete metric space. We denote by $L^2(\Omega, Y)$ the space of all essentially separably valued Borel maps $u\colon \Omega\to Y$ such that for some and thus every $y_0\in Y$ we have $$\int_\Omega d_Y^2(y_0, u(x))\,d\hm^n(x)<\infty.$$ The (Korevaar-Schoen) energy of $u\in L^2(\Omega, Y)$ is defined as follows. 
 For $\varepsilon>0$ we set $$e_\varepsilon(x):= n\cdot \vint_{S(x, \varepsilon)} \frac{d_Y^2(u(x), u(x'))}{\varepsilon^2}\,d\hm^{n-1}(x')$$ whenever $x\in \Omega$ satisfies $d(x, \partial \Omega)> \varepsilon$ and $e_\varepsilon(x)= 0$ otherwise.  The map $u$ is said to belong to $W^{1,2}(\Omega, Y)$ if its energy, defined by
\begin{equation}\label{eq:Korevaar-Schoen-energy}
E(u):= \sup_{f\in C_c(\Omega), \, 0\leq f\leq 1}\left(\limsup_{\varepsilon\to 0} \int_\Omega f(x) e_\varepsilon(x)\,d\hm^n(x)\right),
\end{equation}
 is finite. 
If $u\in W^{1,2}(\Omega, X)$ then there exists a function $e_u\in L^1(\Omega)$, called the energy density function of $u$, such that $e_\varepsilon \,d\hm^n \rightharpoonup e_u\,d\hm^n$ as $\varepsilon\to 0$ and $$E(u) = \int_\Omega e_u(x)\,d\hm^n(x),$$ see \cite[Theorems 1.5.1 and 1.10]{KS93}. In the case that $Y$ is a Riemannian manifold and $u$ is smooth the energy defined in \eqref{eq:Korevaar-Schoen-energy} coincides with the usual energy as defined for example in \cite{BH-2-17}. 

The trace of a Sobolev map $u\in W^{1,2}(\Omega, Y)$ is denoted $\trace(u)$, see \cite[Definition 1.12]{KS93} for the definition. We  mention here that if $u$ has a continuous representative which has a continuous extension to $\overline{\Omega}$, again denoted $u$, then $\trace(u) = u|_{\partial \Omega}$.

\bd\label{def:harmonic}
 A map $u\in W^{1,2}(\Omega, Y)$ is said to be energy minimizing harmonic if $E(u) \leq E(v)$ for all $v\in W^{1,2}(\Omega, Y)$ with $\trace(v) = \trace(u)$. A map $u\colon X\to Y$ is called energy minimizing harmonic if its restriction to every bounded, open ball is energy minimizing harmonic.
\ed

It is well-known that if $Y$ is a Hadamard manifold then a map $u\in W^{1,2}(\Omega, Y)$ is energy minimizing harmonic in the sense above if and only $u$ is a harmonic map in the classical sense (vanishing tension field), see \cite{SU82} and \cite{Sch84}.

Now, let $Y$ be a ${\rm CAT}(0)$-space. It follows from \cite[Theorem 2.2]{KS93} that for every Lipschitz map $f\colon \partial \Omega\to Y$ there exists a unique energy minimizing harmonic map $u\in W^{1,2}(\Omega, Y)$ with $\trace(u) = f$. 
By \cite[Theorem 2.4.6]{KS93} and \cite{Ser94}, the map $u$ is locally Lipschitz continuous in $\Omega$ and H\"older continuous up to the boundary, in particular $u|_{\partial \Omega} = f$.

If $u\in W^{1,2}(\Omega, Y)$ is energy minimizing harmonic then, by  \cite[Lemma 10.2]{EF01}, for every $y_0\in Y$ the function $h\colon \Omega\to \R$ given by $h(z):= d_Y(y_0, u(z))$ is weakly subharmonic in the sense that $\Delta h \geq 0$ weakly. Recall that a function $h\in W^{1,2}(\Omega)$ is said to satisfy $\Delta h \geq \rho$ weakly for some function $\rho\in L^1(\Omega)$ if $$- \int_\Omega \langle \nabla h, \nabla \varphi\rangle\,d\hm^n \geq \int_\Omega  \rho \varphi\,d\hm^n$$ for all non-negative $\varphi\in C^\infty_0(\Omega)$. By \cite[Theorem 1]{Her64}, a continuous and weakly subharmonic function $h\colon \overline{\Omega}\to\R$ with $h|_{\partial \Omega}\leq 0$ satisfies $h\leq 0$ on $\overline{\Omega}$.

The following result will be used in the proof of Theorem~\ref{thm:main-intro}.

\bp\label{prop:interior-grad-estimate-image-ball}
 Let $X$ be a Hadamard manifold with sectional curvature bounded from below and let $B= B(x,r)$ be an open ball in $X$. Let $u\colon B\to Y$ be an energy minimizing harmonic map into some ${\rm CAT}(0)$-space $Y$. If the image of $u$ lies in some ball of radius $R$ then $u$ is $CR$-Lipschitz on the ball $\bar{B}(x, r/3)$, where $C\geq 1$ only depends on $r$, the lower bound on sectional curvature of $X$, and the dimension of $X$.
\ep

\begin{proof}
 Suppose the curvature of $X$ is bounded by $-b^2\leq K_X\leq 0$ for some $b>0$. By \cite[Theorem 1.4]{ZZZ17}, there exists a constant $C_1$ depending only on $r$, the dimension $n$ of $X$, and $b$ such that $u$ is $\lambda$-Lipschitz on $\bar{B}(x, r/3)$ with $\lambda\leq C_1\cdot E(u|_{B(x,s)})^{\frac{1}{2}}$, where we have set $s:= \frac{2r}{3}$. 
 
 It thus suffices to show that $E(u|_{B(x,s)})$ is bounded by $R^2$ times a constant depending only on $r$, $n$, and $b$. Let $y\in Y$ be such that the image of $u$ lies in the ball $\bar{B}(y, R)$. It is not difficult to show that there exists a smooth function $\eta\colon X\to \R$ supported in $B(x,r)$ with $0\leq \eta \leq 1$ everywhere, such that $\eta=1$ on $B(x,s)$, and $|\Delta \eta|\leq K$ everywhere for some constant $K$ depending only on $r$, $b$, and $n$.
 
By \cite[Equation (6)]{Ser94}, we have $\Delta d_Y^2(y, u(x')) \geq 2e_u(x')$ weakly, where $e_u$ denotes the energy density of $u$. We thus obtain
\begin{equation*}
 \begin{split}
  2 E(u|_{B(x,s)}) &\leq 2\int_{B(x,r)} \eta(x') e_u(x')\,d\hm^n(x')\\
   &\leq \int_{B(x,r)}\Delta\eta(x')\cdot d_Y^2(y, u(x'))\,d\hm^n(x')\\
   &\leq K\cdot R^2\cdot \hm^n(B(x, r)).
 \end{split}
\end{equation*}
It follows that the Lipschitz constant $\lambda$ of $u$ on $\bar{B}(x, r/3)$ is bounded by $$\lambda\leq C_1\cdot E(u|_{B(x,s)})^{\frac{1}{2}} \leq C R$$ for some constant $C$ depending on $r$, $b$, and $n$. This completes the proof.
\end{proof}

\section{Lipschitz quasi-isometric maps}\label{sec:Lip-qi-map}

We will need:

\bp\label{prop:Lipschitz-quasi-isometric}
 Let $X$ be a Hadamard manifold with sectional curvature bounded from below and let $Y$ be a ${\rm CAT}(0)$-space. Then every quasi-isometric map $f\colon X\to Y$ is at finite distance from a quasi-isometric map $\tilde{f}\colon X\to Y$ which is moreover Lipschitz.
\ep


We first show the following lemma which will also be used later.

\bl\label{lem:Hadamard-bdd-growth}
 Let $X$ be a Hadamard manifold with sectional curvature bounded from below. Then for every $0<r<R<\infty$ there exists $N\in\N$ such that every ball in $X$ of radius $R$ can be covered by $N$ balls of radius $r$. 
\el

\begin{proof}
 Fix $0<r<R<\infty$ and $x\in X$. Let $A\subset \bar{B}(x,R)$ be a maximally $r$-separated subset. Thus, distinct points in $A$ have distance at least $r$ and the union of open $r$-balls centered at points in $A$ covers the ball $\bar{B}(x, R)$. The open balls centered at points in $A$ and with radius $\frac{r}{2}$ are pairwise disjoint and contained in the ball $B(x, R+\frac{r}{2})$. Let $a_1, \dots, a_k\in A$ be distinct points and denote by $m$ the volume of the ball of radius $\frac{r}{2}$ in Euclidean $\R^n$, where $n$ is the dimension of $X$. Volume comparison with Euclidean space yields $\hm^n(B(a_i, \frac{r}{2}))\geq m$ for every $i$, see \cite[Theorem 3.101]{GHL04}. By the same theorem, applied to a model space of constant negative curvature, we obtain $$\hm^n\left(B\left(x, R+\frac{r}{2}\right)\right)\leq M$$ for some $M$ depending only on $R$, $n$, and the lower bound on sectional curvature. We conclude that $$k\cdot m \leq \sum_{i=1}^k \hm^n\left(B\left(a_i, \frac{r}{2}\right)\right) = \hm^n\left(\bigcup_{i=1}^k B\left(a_i, \frac{r}{2}\right)\right) \leq \hm^n\left(B\left(x, R+\frac{r}{2}\right)\right)\leq M$$ and hence that $k\leq \frac{M}{m}$. This shows that $A$ has at most $\frac{M}{m}$ points. Since the union of the open $r$-balls centered at points in $A$ covers the ball $\bar{B}(x, R)$ the proof is complete.
\end{proof}

We now prove Proposition~\ref{lem:Hadamard-bdd-growth}

\begin{proof}
 Let $f\colon X\to Y$ be an $(L, c)$-quasi-isometric map and let $Z\subset X$ be a maximally $1$-separated subset of $X$. 
 It is easy to show that the family of balls given by $\left\{\bar{B}(z, 4): z\in Z\right\}$ has bounded multiplicity. Indeed, let $x\in X$ and let $z_1,\dots, z_k\in Z$ be distinct points such that $d(x,z_i)\leq 4$ for all $i$. By Lemma~\ref{lem:Hadamard-bdd-growth} the ball $\bar{B}(x, 4)$ can be covered by $N$ open balls of radius $\frac{1}{4}$, where $N$ only depends on the lower bound on sectional curvature and the dimension of $X$. Since each of these balls can contain at most one element of $Z$ our claim follows.
 
Now, the restriction $f|_Z$ of $f$ to $Z$ is $(L+c)$-Lipschitz. Moreover, $Y$ is Lipschitz $k$-connected for every $k\in\N$. Thus, \cite[Lemma 5.3]{Wen06} implies that the map $f|_Z$ has a Lipschitz extension $\tilde{f}\colon X\to Y$ whose Lipschitz constant only depends on $N$. By the triangle inequality, the map $\tilde{f}$ is at bounded distance from $f$ and hence also quasi-isometric. This concludes the proof.
\end{proof}

\section{The boundary estimate}\label{sec:boundary-estimate}

Let $(X, d_X)$ be a Hadamard manifold of dimension $n\geq 2$ and of pinched negative curvature $-b^2 \leq K_X\leq -a^2$ for some $a,b>0$. Let $(Y, d_Y)$ be a ${\rm CAT}(0)$-space which is locally compact and Gromov hyperbolic. Suppose $f\colon X\to Y$ is a quasi-isometric map which is moreover Lipschitz. Thus there exist $L\geq 1$ and $c>0$ such that
\begin{equation}\label{eq:f-qi-Lipschitz}
 L^{-1}\cdot d_X(x, x') - c \leq d_Y(f(x), f(x')) \leq L\cdot d_X(x, x')
\end{equation}
for all $x, x'\in X$.
Let $x_0\in X$ and set $\bar{B}_R:= \bar{B}(x_0, R)$ whenever $R>0$. We furthermore set $S_R:= S(x_0, R)$. There exists a unique continuous energy minimizing harmonic map $u_R\colon \bar{B}_R \to Y$ which coincides with $f$ on $S_R$, see Section~\ref{sec:Sobolev}. The main aim of this section is to establish:

\bp\label{prop:boundary-estimate}
 There exist constants $\alpha, \beta\geq 1$ such that for every $R>0$ and $x\in \bar{B}_R$ we have $$d_Y(f(x), u_R(x)) \leq \alpha \cdot d_X(x, S_R) + \beta.$$
\ep

The proof of the analogous result \cite[Proposition 3.7]{BH-2-17} when $Y$ is a Hadamard manifold with curvature bounded from below heavily depends on the existence, established in \cite[Proposition 2.4]{BH-2-17}, of a smooth map at finite distance of $f$ with bounded first and second covariant derivative. In the singular setting we work in, such a result is of course not available. We circumvent this problem by using the following lemma.

\bl\label{lem:Bonk-Schramm-for-fX}
 The set $f(X)$, equipped with the metric from $Y$, admits a rough-isometric map $\psi\colon f(X)\to\Hyp^k$ for some $k\in\N$.
\el

Recall that a map $\psi\colon Z\to W$ between metric spaces $(Z, d_Z)$ and $(W, d_W)$ is called $(\lambda, c)$-rough-isometric if $$\lambda \cdot d_Z(z,z') - c \leq d_W(\psi(z), \psi(z')) \leq \lambda\cdot d_Z(z,z') + c$$ for all $z,z'\in Z$.
The idea is to use the well-known Bonk-Schramm embedding theorem \cite{BS00}. We cannot use their embedding theorem directly since $f(X)$ is not geodesic. We will therefore use injective hulls of metric spaces. See Section~\ref{sec:injective-hull} above and \cite{Lan13} for the definition as well as the properties we need. 

\begin{proof}
 Denote by $Z$ the set $f(X)$ equipped with the induced metric from $Y$. Denote by $E(Z)$ the injective hull of $Z$ and recall that $E(Z)$ is a  geodesic metric space and that $Z$ embeds isometrically into $E(Z)$. Since $Z$ is Gromov hyperbolic (as a subset of $Y$) it follows from \cite[Proposition 1.3]{Lan13} that $E(Z)$ is also Gromov hyperbolic. 
 
 We claim that $E(Z)$ is in finite Hausdorff distance of $Z$. For this, notice first that the space $E(Z)$ can be viewed as a subset of the injective hull $E(Y)$ of $Y$, see Section~\ref{sec:injective-hull}. Now, let $p$ be a point in $E(Z)\subset E(Y)$. The proof of \cite[Proposition 1.3]{Lan13} shows that $p$ lies at distance at most $\delta_1$ from a geodesic $[z,z']$ in $Y$ between two points $z=f(x)$ and $z'=f(x')$ in $Z$, where $\delta_1$ only depends on the Gromov hyperbolicity constant of $Y$. Let $[x,x']$ be the geodesic segment in $X$ from $x$ to $x'$. By the stability of quasi-geodesics \cite[Theorem III.H.1.7]{BH99}, the quasi-geodesic $f([x,x'])$ is at distance at most $\delta_2$ from the geodesic $[z,z']$, where $\delta_2$ only depends on the Gromov hyperbolicity constant of $Y$ and the quasi-isometric constants of $f$. Thus $p$ lies at distance at most $\delta_1 + \delta_2$ from a point in $Z$. This proves our claim.
 
Since $E(Z)$ is at finite distance from $Z$ and $f$ is quasi-isometric it easily follows from Lemma~\ref{lem:Hadamard-bdd-growth} that $E(Z)$ has bounded growth at some scale as defined in \cite{BS00}. That is, there exist $0<r<R<\infty$ and $N\in\N$ such that every ball of radius $R$ in $E(Z)$ can be covered by at most $N$ balls of radius $r$. Since $E(Z)$ is also geodesic and Gromov hyperbolic it follows from the Bonk-Schramm embedding theorem \cite{BS00} that $E(Z)$ admits a rough-isometric map $\psi\colon E(Z)\to \Hyp^k$ for some $k\in\N$. Since $E(Z)$ contains $Z$ isometrically, the proof is complete.
\end{proof}

We are now ready for the proof of Proposition~\ref{prop:boundary-estimate}. It uses Lemma~\ref{lem:Bonk-Schramm-for-fX} but is otherwise very similar to that of \cite[Proposition 3.7]{BH-2-17}.

\begin{proof}
 Denote by $Z$ the set $f(X)$ equipped with the induced metric from $Y$. By Lemma~\ref{lem:Bonk-Schramm-for-fX} there exists a $(\lambda, \bar{c})$-rough-isometric map $\psi\colon Z\to \Hyp^k$ for some $\lambda, \bar{c}>0$ and $k\in \N$. Since the composition $\psi\circ f$ is quasi-isometric there exist by \cite[Proposition 2.4]{BH-2-17} constants $A$ and $M$ and a smooth map $\tilde{f}\colon X\to \Hyp^k$ such that $$d_{\Hyp^k}(\psi\circ f(x), \tilde{f}(x))\leq M$$ for all $x\in X$ and such that the $\|D\tilde{f}\|\leq A$ and $\|\tau(\tilde{f})\|\leq A^2$. 
  
Fix $R>0$ and $x\in B_R$ and define two continuous functions $\varphi_1, \varphi_2\colon \bar{B}_R\to \R$ by $$\varphi_1(z):= \lambda\cdot d_Y(f(x), u_R(z))$$ and $$\varphi_2(z):= \frac{2nA^2}{a (n-1)}\cdot (d_X(x_0, z)- R).$$ By \cite[Lemma 10.2]{EF01} the function $\varphi_1$ is weakly subharmonic, see also Section~\ref{sec:Sobolev} above. Furthermore, the function $d_{x_0}(z):= d_X(x_0, z)$ satisfies $\Delta d_{x_0}\geq a\cdot (n-1)$ away from $x_0$, see Section~\ref{sec:Riem-prelim} above, and so the function $\varphi_2$ satisfies $\Delta\varphi_2\geq 2nA^2$ weakly.

Now, we define a third function $\varphi_3\colon \bar{B}_R\to \R$ as follows. Set $y_0: = \psi(f(x))$ and embed $\Hyp^k$ isometrically into $\Hyp^{k+1}$. We pick a point $y_1$ on the geodesic in $\Hyp^{k+1}$ passing through $y_0$ perpendicular to $\Hyp^k$ which is sufficiently far from $y_0$ and define $$\varphi_3(z):= - d_{\Hyp^{k+1}}(y_1, \tilde{f}(z)) + d_{\Hyp^{k+1}}(y_0, y_1).$$
From \eqref{eq:hessian-distance-hyp} and \eqref{eq:laplace-composition} we see that the function $\varphi_3$ satisfies 
\begin{equation*}
 |\Delta \varphi_3| = |\Delta(d_{y_1}\circ\tilde{f})| \leq  \|Dd_{y_1}\| \cdot \|\tau(\tilde{f})\| + n \cdot\coth(d_{\Hyp^{k+1}}(y_0, y_1)) \cdot \|D\tilde{f}\|^2
\end{equation*} 
everywhere on $B_R$, where $d_{y_1}\colon\Hyp^{k+1}\to\R$ is given by $d_{y_1}(w):= d_{\Hyp^{k+1}}(y_1, w)$. If $y_1$ is chosen sufficiently far from $y_0$ then it follows that $|\Delta\varphi_3|\leq 2nA^2$ everywhere on $B_R$. Consequently, the continuous function $\varphi\colon \bar{B}_R\to \R$ defined by $\varphi:= \varphi_1+\varphi_2+\varphi_3$ is weakly subharmonic.

We now estimate $\varphi$ on $S_R$. For this let $z\in S_R$ and notice that $\varphi_1(z)= \lambda\cdot d_Y(f(x), f(z))$ and $\varphi_2(z)=0$. Since $y_1\in\Hyp^{k+1}$ is on the geodesic from $y_0$ perpendicular to $\Hyp^k$ it follows from the hyperbolic law of cosines that $$d_{\Hyp^{k+1}}(y, y_1) \geq d_{\Hyp^{k+1}}(y, y_0) + d_{\Hyp^{k+1}}(y_0, y_1) - \log(4)$$ for every $y\in \Hyp^k\subset \Hyp^{k+1}$. From this we conclude that 
\begin{equation*}
 \begin{split}
  \varphi_3(z) &= - d_{\Hyp^{k+1}}(y_1, \tilde{f}(z)) + d_{\Hyp^{k+1}}(y_0, y_1)\\
  &\leq - d_{\Hyp^k}(y_0, \tilde{f}(z)) + \log(4)\\
  &\leq - d_{\Hyp^k}(\psi(f(x)),\psi(f(z))) + M + \log(4)\\
  &\leq - \lambda d_Y(f(x), f(z)) + M',
 \end{split}
\end{equation*}
where $M':= \bar{c} + M + \log(4)$. It follows that $\varphi(z) \leq M'$ for every $z\in S_R$. Since $\varphi$ is weakly subharmonic and continuous we thus obtain from \cite[Theorem 1]{Her64} or from Section~\ref{sec:Sobolev} above that $\varphi(z) \leq M'$ for all $z\in\bar{B}_R$ and, in particular, also for $z=x$. Since $|\varphi_3(x)| \leq M$ we conclude that $$d_Y(f(x), u_R(x)) \leq \frac{2nA^2}{\lambda a (n-1)}\cdot d_X(x, S_R) + \frac{M' +M}{\lambda},$$ which completes the proof.
\end{proof}

\section{Distance between harmonic and quasi-isometric maps}\label{sec:dist-qi-harmonic-unif}

The proof of the following proposition is almost identical to that of \cite[Proposition 3.5]{BH-2-17} except that we use the Gromov product instead of angle estimates. The latter are not available in our setting. Let $(X, d_X)$, $(Y, d_Y)$, $f$, $x_0$, $\bar{B}_R$, and $u_R$ be as in Section~\ref{sec:boundary-estimate}.

\bp\label{prop:max-dist-qi-harmonic}
 There exists $\rho\geq 1$ such that for every $R\geq 1$ we have $$d_Y(u_R(x), f(x)) \leq \rho$$ for all $x\in \bar{B}_R$.
\ep

We turn to the proof and let $a,b>0$ be such that the sectional curvature of $X$ satisfies $-b^2\leq K_X\leq -a^2$. Let $C\geq 1$ be as in Proposition~\ref{prop:interior-grad-estimate-image-ball} for the radius $r=3$. Let $\delta>0$ be such that $Y$ is $\delta$-hyperbolic in the sense of Definition~\ref{def:Gromov-hyp} and let $\delta'>0$ be the constant from Lemma~\ref{lem:div-geodesics}. Denote by $L$ and $c$ the constants from \eqref{eq:f-qi-Lipschitz} and by $c'$ the constant from Lemma~\ref{lem:effecxt-qi-Gromov-product}. Let $\alpha, \beta\geq 1$ be as in Proposition~\ref{prop:boundary-estimate}. Finally, let $M$ and $N$ be the constants appearing in the uniform estimates on the harmonic measure on distance spheres in $X$ proved in \cite[Theorem 1.1]{BH-harmonic-measure-18}.

We choose $T>3$ so large that inequality \eqref{eq:choice-of-T} below holds and that $\gamma$, as defined in \eqref{eq:def-angle-gamma} below, satisfies $\gamma<\frac{\pi}{2}$.
We argue by contradiction and assume Proposition~\ref{prop:max-dist-qi-harmonic} is false. There then exists a sequence $R_k\to\infty$ such that 
\begin{equation}\label{eq:distance-large-rho}
 \rho_k:= \sup\left\{x\in \bar{B}_{R_k}: d_Y(u_{R_k}(x), f(x))\right\} \to \infty
\end{equation}
as $k\to\infty$. We now abbreviate $u_k:= u_{R_k}$. Let $k\geq 1$ be sufficiently large so that $$\rho_k> \max\left\{2T\alpha + \beta, 2LT + 6\delta', 4LMT\gamma^{-N}\right\}$$ and so that $\rho_k$ satisfies inequality \eqref{eq:rhok-divergence-choice} below. Since $u_k$ and $f$ are continuous on $\bar{B}_{R_k}$ the supremum in \eqref{eq:distance-large-rho} is achieved at some point $x\in \bar{B}_{R_k}$. By Proposition~\ref{prop:boundary-estimate} and the choice of $\rho_k$ we have $$d_X(x, S_{R_k}) \geq \frac{\rho_k - \beta}{\alpha}> 2T.$$ In particular, the ball $\bar{B}(x, 2T)$ is contained in $B_{R_k}$. We first prove:

\bl\label{lem:harm-Lip-annulus}
The map $u_k$ is $2C\rho_k$-Lipschitz on $\bar{B}(x, T)$ and satisfies
\begin{equation}\label{eq:fx-ukz-bounds-and-Lip}
\frac{\rho_k}{2} \leq d_Y(f(x), u_k(z)) \leq \rho_k + LT
\end{equation}
 for all $z\in \bar{B}(x, T)$.
\el

\begin{proof}
 For every $z\in X$ with $d(x,z)\leq 2T$ we have $$d_Y(f(x), u_k(z)) \leq d_Y(f(x), f(z)) + d_Y(f(z), u_k(z)) \leq Ld(x,z) + \rho_k,$$ which implies in particular the second inequality in \eqref{eq:fx-ukz-bounds-and-Lip} and that $u(\bar{B}(x, 2T))\subset \bar{B}(f(x), 2\rho_k)$ because $2LT< \rho_k$. Now, let $z\in \bar{B}(x, T)$. Since $\bar{B}(z, 3)\subset \bar{B}(x, 2T)$ it thus follows from Proposition~\ref{prop:interior-grad-estimate-image-ball}, applied with $r=3$, that $u_k$ is $2C\rho_k$-Lipschitz on the ball $\bar{B}(z,1)$ and hence also on the ball $\bar{B}(x, T)$ since balls in $X$ are geodesic.

It remains to verify the first inequality in \eqref{eq:fx-ukz-bounds-and-Lip}. Suppose it does not hold everywhere. Then there exists $z_1\in \bar{B}(x, T)$ such that $$h(z_1):= d_Y(f(x), u_k(z_1)) = \frac{\rho_k}{2}.$$ Set $r_1:= d(x, z_1)>0$. The Lipschitz continuity just proved implies $$h(z) \leq h(z_1) + d_Y(u_k(z_1), u_k(z)) \leq \frac{3\rho_k}{4}$$ for all $z$ in the set $\Sigma:=S(x, r_1)\cap \bar{B}(z_1, \frac{1}{8C})$. Using the hyperbolic law of cosines and comparing with the hyperbolic plane of curvature $-b^2$ we see that $\Sigma$ contains the intersection of $S(x, r_1)$ with a geodesic cone $C_\gamma$ based at $x$ and with angle 
\begin{equation}\label{eq:def-angle-gamma}
 \gamma= \frac{\sqrt{\cosh(\frac{b}{8C}) - 1}}{\sinh(b T)}.
\end{equation}
Let $\sigma$ denote the harmonic measure on $S(x, r_1)$. See \cite{BH-harmonic-measure-18} for the definition. Since $h$ is continuous and weakly subharmonic the harmonic function $\xi$ on $\bar{B}(x, r_1)$ which equals $h$ on $S(x, r_1)$ satisfies $$\rho_k = h(x) \leq \xi(x) = \int_{S(x, r_1)}\xi\,d\sigma = \int_{S(x, r_1)}h\,d\sigma$$ and hence $$\int_{S(x, r_1)}(h-\rho_k)\,d\sigma\geq 0.$$
Since $h-\rho_k\leq LT$ on $S(x, r_1)$  and $h-\rho_k \leq -\frac{\rho_k}{4}$ on $C_\gamma\cap S(x,r_1)$ it follows that $\sigma(C_\gamma)\leq \frac{4LT}{\rho_k}$. From the uniform lower bound on the harmonic measure of geodesic cones proved in \cite[Theorem 1.1]{BH-harmonic-measure-18} we thus obtain $$\frac{1}{M}\cdot \gamma^N \leq \sigma(C_\gamma)\leq \frac{4LT}{\rho_k},$$ which contradicts the choice of $\rho_k$. The proof is complete.
\end{proof}

We now define a subset $U\subset X$ by $$U:= \left\{z\in S(x, T): d_Y(f(x), u_k(z))\geq \rho_k - \frac{T}{2L}\right\}$$ and prove:

\bl\label{lem:lower-bound-Gromov-product-f}
 For all $z_1, z_2\in U$ we have $$(f(z_1) \mid f(z_2))_{f(x)} \geq \frac{T}{4L} - \frac{c}{2} - 2\delta.$$
\el

\begin{proof}
  Let $z_1, z_2\in U$ and notice that for $i=1, 2$, we have
  \begin{equation*}
   \begin{split} 
    2\cdot (f(z_i)\mid u_k(z_i))_{f(x)} &= d_Y(f(x), f(z_i)) + d_Y(f(x), u_k(z_i)) - d_Y(f(z_i), u_k(z_i))\\
     &\geq \frac{T}{L} - c + \rho_k - \frac{T}{2L}  - \rho_k\\
     &= \frac{T}{2L} - c.
   \end{split}
  \end{equation*}
We next claim that 
\begin{equation}\label{eq:lower-bound-Gromov-prod-uk}
 (u_k(z_i)\mid u_k(x))_{f(x)} \geq \frac{\rho_k}{4} - \frac{3\delta'}{2}.
\end{equation}
In order to show this, fix $i$ and set $y:= f(x)$, $y_1:= u_k(z_i)$, and $y_2:= u_k(x)$ and recall that $d_Y(y, y_1) \geq \rho_k- \frac{T}{2L}\geq \frac{\rho_k}{2}$ and $d_Y(y, y_2) = \rho_k$. For $j=1,2$, let $y'_j$ be the point on the geodesic from $y$ to $y_j$ with $d_Y(y, y'_j) = \frac{\rho_k}{4}$. Let $\xi$ be the geodesic in $X$ from $x$ to $z_i$. By Lemma~\ref{lem:harm-Lip-annulus}, the curve $u_k\circ \xi$ stays outside the ball $B(y, \frac{\rho_k}{2})$ and has length bounded from above by $2CT\rho_k$. Since $\rho_k$ was chosen so large that
\begin{equation}\label{eq:rhok-divergence-choice}
 \rho_k + LT + 2CT\rho_k < 2^{\frac{\rho_k - 4}{4\delta'}}
\end{equation}
it follows from Lemma~\ref{lem:div-geodesics} that $d_Y(y'_1, y'_2) \leq 3\delta'$. This is easily seen to imply \eqref{eq:lower-bound-Gromov-prod-uk}, which proves our claim.

Finally, we use the definition of $\delta$-hyperbolicity of $Y$, the estimates above, and the fact that $\frac{\rho_k}{4} - \frac{3\delta'}{2}\geq \frac{T}{4L} - \frac{c}{2}$ to conclude that 
\begin{equation*}
  (f(z_i)\mid u_k(x))_{f(x)} \geq \min\left\{(f(z_i)\mid u_k(z_i))_{f(x)}, (u_k(z_i)\mid u_k(x))_{f(x)}\right\} - \delta \geq \frac{T}{4L} - \frac{c}{2} - \delta
 \end{equation*}
and hence
\begin{equation*}
  (f(z_1)\mid f(z_2))_{f(x)} \geq \min\left\{(f(z_1)\mid u_k(x))_{f(x)}, (f(z_2)\mid u_k(x))_{f(x)}\right\} - \delta \geq \frac{T}{4L} - \frac{c}{2} - 2\delta.
 \end{equation*}
 This completes the proof.
\end{proof}

The next lemma provides a contradiction to the previous lemma since we had chosen $T$ so large that 
\begin{equation}\label{eq:choice-of-T}
\frac{L}{a}\cdot \log\left(4 M^N(2L^2+1)^N \right) + c'< \frac{T}{4L} - \frac{c}{2} - 2\delta.
\end{equation}
The lemma will thus finish the proof of Proposition~\ref{prop:max-dist-qi-harmonic}.

\bl\label{lem:upper-bound-Gromov-prod}
 There exist $z_1, z_2\in U$ such that $$(f(z_1) \mid f(z_2))_{f(x)} \leq \frac{L}{a}\cdot \log\left(4 M^N(2L^2+1)^N \right) + c'.$$
\el

\begin{proof}
  Denote by $\sigma$ the harmonic measure on $S(x,T)$. Let $h$ be the continuous and weakly subharmonic function given by $h(z):= d_Y(f(x), u_k(z))$. Comparing with a harmonic function exactly as in the proof of Lemma~\ref{lem:harm-Lip-annulus} we obtain that $$\int_{S(x,T)} (h-\rho_k)\,d\sigma \geq 0.$$
By the definition of $U$ and by Lemma~\ref{lem:harm-Lip-annulus}, we have $h(z)- \rho_k \leq LT$ for all $z\in S(x, T)$ and $h(z) - \rho_k < - \frac{T}{2L}$ whenever $z\in S(x, T) \setminus U$. This together with the above integral inequality yields $$\sigma(U) \geq \frac{1}{2L^2+1}.$$
The uniform upper bound on the harmonic measure proved in \cite[Theorem 1.1]{BH-harmonic-measure-18} thus shows that there exist $z_1, z_2\in U$ such that the angle $\gamma'$ between them, as seen from the point $x$, satisfies $$\gamma'\geq \left(\frac{\sigma(U)}{M}\right)^N \geq \frac{1}{M^N\cdot (2L^2+1)^N}.$$ From this, \cite[Lemma 2.1]{BH-1-17}, and Lemma~\ref{lem:effecxt-qi-Gromov-product} it follows that $$(f(z_1)\mid f(z_2))_{f(x)} \leq L (z_1\mid z_2)_x + c' \leq  \frac{L}{a}\cdot \log\left(4 M^N(2L^2+1)^N \right) + c',$$ which concludes the proof.
\end{proof}

\section{Completing the proof of the main theorem}

We complete the proof of Theorem~\ref{thm:main-intro}. 
Let $(X, d_X)$ and $(Y, d_Y)$ be spaces as in the statement of the theorem and let $f\colon X\to Y$ be a quasi-isometric map. By Proposition~\ref{prop:Lipschitz-quasi-isometric} we may assume that $f$ is also $L$-Lipschitz continuous for some $L>0$. Fix a basepoint $x_0\in X$ and set $B_R:= B(x_0, R)$ and $S_R:= S(x_0, R)$ whenever $R>0$. Let furthermore $u_R\colon \bar{B}_R\to Y$ be the unique continuous energy minimizing harmonic map which coincides with $f$ on $S_R$, see Section~\ref{sec:Sobolev}. Proposition~\ref{prop:max-dist-qi-harmonic} shows that there exists $\rho$ such that
\begin{equation}\label{eq:dist-qi-harmonic-R}
 d_Y(u_R(x), f(x)) \leq \rho
\end{equation}
for all $R\geq 1$ and every $x\in \bar{B}_R$. From this and the Lipschitz continuity of $f$ it follows that for every $x\in B(x_0, R-4)$ the image of $u_R(B(x, 3))$ is contained in a ball of radius $\rho + 3L$. Proposition~\ref{prop:interior-grad-estimate-image-ball} implies that $u_R$ is $L'$-Lipschitz on $B(x,1)$ for some $L'$ which does not depend on $x$ or $R$. Consequently, $u_R$ is $L'$-Lipschitz on $B(x_0, R-4)$. 

Fix a sequence $R_k\to\infty$ and set $u_k:= u_{R_k}$. By Arzela-Ascoli theorem, a diagonal subsequence argument, and by \eqref{eq:dist-qi-harmonic-R} we may thus assume that there exists an $L'$-Lipschitz map $u\colon X\to Y$ such that $u_k$ converges to $u$ uniformly on compact sets and that $ d_Y(u(x), f(x)) \leq \rho$ holds for every $x\in X$. 

It remains to show that $u$ is energy minimizing harmonic. Fix $s>0$. The restriction of $u$ to $B_s$ is in $W^{1,2}(B_s, Y)$ since $u$ is Lipschitz. Now, suppose there exist $\varepsilon>0$ and $v\in W^{1,2}(B_s, Y)$ such that $\trace(v) = u|_{S_s}$ and $E(v)\leq E(u|_{B_s}) - \varepsilon$. Let $\delta\in(0,1)$ be sufficiently small, to be determined later. For $k$ sufficiently large (depending on $\delta$) the map $h\colon S_s \cup S_{s+\delta}\to Y$ defined by $h= u$ on $S_s$ and $h= u_k$ on $S_{s+\delta}$ is $2L'$-Lipschitz. Since the ball $B_{s+1}$ is doubling and hence also $A_\delta:= \bar{B}_{s+\delta}\setminus B_s$ and $Y$ is Lipschitz $m$-connected for every $m$ it follows from \cite[Theorem 1.5]{LS05} that $h$ has an $L''$-Lipschitz extension $\bar{h}\colon A_\delta\to Y$ with Lipschitz constant $L''$ not depending on $\delta$ or $k$. We now define a map $v_k\colon \bar{B}_{s+\delta}\to Y$ as follows. For $x\in B_s$ set $v_k(x):= v(x)$ and for $x\in A_\delta$ set $v_k(x) = \bar{h}(x)$. Then $v_k\in W^{1,2}(B_{s+\delta}, Y)$ with $\trace(v_k) = u_k|_{S_{s+\delta}}$, see \cite[Theorem 1.12.3]{KS93}. Since $u_k$ is energy minimizing harmonic we have
\begin{equation*}
 \begin{split}
E((u_k)|_{B_{s+\delta}})&\leq E(v_k)\leq E(v) + n (L'')^2\cdot \hm^n(A_\delta)\\
 &\leq  E(u|_{B_{s+\delta}}) - \varepsilon + n(L'')^2\cdot \hm^n(A_\delta).
 \end{split}
\end{equation*}
However, the right-hand side is strictly smaller than $E(u|_{B_{s+\delta}}) - \frac{\varepsilon}{2}$ whenever $\delta>0$ is sufficiently small and $k$ is sufficiently large. This contradicts the lower semi-continuity of the energy \cite[Theorem 1.6.1]{KS93}. We conclude that $u$ is indeed energy minimizing harmonic. This completes the proof.

\end{document}